\documentclass[12pt]{article}
\newtheorem{thm}{Theorem}

\newtheorem{example}{Example}

\font\msbm=msbm10 scaled 1200
\font\msbmscript=msbm8
\font\msbmscriptscript=msbm6
\newfam\Bbbfam
\textfont\Bbbfam=\msbm
\scriptfont\Bbbfam=\msbmscript
\scriptscriptfont\Bbbfam=\msbmscriptscript
\def\Bbb#1{{\fam\Bbbfam#1}}

\newcommand{\mod}{\mathop{\rm mod}\nolimits}

\makeatletter
\@addtoreset{equation}{section}
\def\thethm{\arabic{thm}\protect\@blinkpoint}
\def\thecorollary{\arabic{corollary}\protect\@blinkpoint}
\def\theexample{\arabic{example}\protect\@blinkpoint}
\def\thesection{\arabic{section}\protect\@blinkpoint}

\def\@blinkpoint{.}
\let\@blinkref=\ref
\def\ref#1{{\def\@blinkpoint{}\@blinkref{#1}}}
\let\@afterindentfalse\@afterindenttrue
\makeatother

\def\halmos{\hfill\rule{6pt}{6pt}}
\begin{document}
\thispagestyle{empty}
\medskip
\begin{center}
  {\large\bf Not all limit points of poles of the Pad\'e
approximants are obstructions for poinwise convergence}
\footnote{This work was supported by Russian Foundation for Basic
Research under Grant \# 04-01-96006.}

\medskip
Victor M. Adukov \\
{ \it Department of Differential Equations and Dynamical Systems,\\
Southern Ural State University, Lenin avenue 76, 454080 Chelyabinsk, Russia\\
E-mail address: avm@susu.ac.ru}
\end{center}
\medskip
{\small \it
\noindent{\bf Abstract}

In the work it is shown that not all limit points of poles of the
Pad\'e approximants for the last intermediate row are obstructions
for pointwise convergence of the whole row to an approximable
function. The corresponding examples are constructed. }
\medskip

The notion of intermediate rows was introduced in the
work~\cite{Sidi}, where sufficient conditions for the convergence of
the whole intermediate row were obtained. We will use the explicit
description of the set of limit points of poles of the Pad\'e
approximants for the last intermediate row obtained
in~\cite{Adukov}.

The paper is the continuation of the works ~\cite{Adukov}.

We will construct the examples in the class of rational functions. Let
$r(z)=\frac{N(z)}{D(z)}$ be a  strictly proper rational  fraction. Here
$N(z), D(z)$ are coprime polynomials, $\deg D(z)=\lambda$, and all poles
of $r(z)$  lie in the  disk $|z|<R$.

Let $\left(V_k(z), U_k(z)\right)$ be the unique solution of the Bezout
equation
$$
N(z)V_k(z)+D(z)U_k(z)=z^k,\ \ \ k\ge 0,
$$
where the polynomial $V_k(z)$ satisfy the inequality $\deg V_k(z)<\lambda$.
This solution is called the minimal solution of the Bezout
equation. It is easily seen that for   $n\ge 0$ the polynomials
$V_{n+\lambda}(z), -U_{n+\lambda}(z)$  are the denominator
and numerator of the Pad\'e approximants  $\pi_{n,\lambda-1}^r(z)$ of
type  $(n,\lambda-1)$ for $r(z)$, respectively:
$$
\pi_{n,\lambda-1}^r(z)=-\frac{U_{n+\lambda}(z)}{V_{n+\lambda}(z)}
$$
(see ~\cite{Adukov}, Theorem 4.1).
For $V_k(z)$ there is the following
explicit formula:
$$
V_k(z)=\sum_{i=1}^{\lambda}\sum_{j=0}^{i-1}d_iv_{k+i-j-1}z^j
$$
(see~\cite{Adukov}, formula~(4.8)).
Here $D(z)= z^\lambda + d_{\lambda-1}z^{\lambda-1}+\cdots +d_0$ and
$v_k$ is the coefficient of $z^{\lambda-1}$ in the polynomial $V_k(z)$.

We will assume that {\em all roots $z_1,\ldots,z_\lambda$
of the polynomial $D(z)$  are  simple}. Then $v_k$ is found as follows
$$
v_k=C_1z_1^k+\cdots +C_\lambda z_\lambda^k
$$
(see ~\cite{Adukov}, formula (4.10)), the numbers $C_j$ are defined in
this case by the formula
$$
C_j=\frac{1}{D_j^2(z_j)A_j},
$$
where  $D_j(z)=\frac{D(z)}{z-z_j}$ and  $A_j$ is  a residue  of $r(z)$ at $z_j$
(see~\cite{Adukov}, formula  (4.14)). Substituing $v_k$ into $V_k(z)$,
we obtain
\begin{equation} \label{1}
V_k(z)= \sum_{j=1}^{\lambda}C_jz_j^kD_j(z).
\end{equation}

We will also suppose that {\em the dominant poles
$z_1,\ldots,z_\nu, \ 1<\nu< \lambda-1$,  are vertices of a regular $\sigma$-gon,
$\sigma\ge\nu$, and}
$
\rho\equiv |z_1| =\ldots=|z_{\nu}|>|z_{\nu +1}|>|z_{\nu +2}|\ge
\ldots \ge |z_{\lambda}|.
$

Let us introduce the following family of polynomials:
$$
\omega_m(z)=\sum_{j=1}^{\nu}C_jz_j^m\Delta_k(z),
\ \ \ m=0,1,\ldots,\sigma-1,
$$
where $\Delta(z)=(z-z_1)\cdots (z-z_\nu)$, $\Delta_k(z)=\frac{\Delta(z)}{z-z_k}$.

Let $n+\lambda\equiv m(\mod \sigma)$. If we divide  the sum in (\ref{1}) by
two sums over the dominant and nondominant poles, we get
\begin{equation} \label{2}
V_{n+\lambda}(z)= (z-z_{\nu+1})\cdots (z-z_\lambda)z_1^{n+\lambda-m}\omega_m(z)
+D(z)\sum_{j=\nu+1}^\lambda\frac{C_jz_j^{n+\lambda}}{z-z_j}.
\end{equation}
Here we take into account the equality $z_1^{n+\lambda-m}=\cdots =
z_\nu^{n+\lambda-m}$ for $n+\lambda\equiv m(\mod \sigma)$.

From the representation~(\ref{2}) we see that for $n\to \infty$,
$n\in \Lambda_m\equiv
\left\{ n\in \Bbb N \Bigl| n+\lambda\equiv m(\mod \sigma)\right\}$
there exists
$$
\lim z_1^{-(n+\lambda-m)}V_{n+\lambda}(z)=
(z-z_{\nu+1})\cdots (z-z_\lambda)\omega_m(z).
$$
This means that the set of limit points of poles of the Pad\'e approximants
$\left\{\pi^r_{n,\lambda-1}(z)\right\}_{n\in\Lambda_m}$ consists of the
poles $z_{\nu+1},\ldots,z_\lambda$ and the additional limit points that
are the roots of the polynomials $\omega_m(z), \ \ m=0,1,\ldots,\sigma-1$
(see also Theorem 2.7 in~\cite{Adukov}). The additional limit points are
different from the dominant poles
$z_1,\ldots,z_\nu$ (see Proposition 2.1 in~\cite{Adukov}), but they can
be coincided with the poles $z_{\nu+1},\allowbreak \ldots,z_\lambda$.
The corresponding example  will be given below. It should be noted that
every additional limit point $\zeta_0$ generates the sequence of defects
of the Pad\'e approximants $\pi^r_{n,\lambda-1}(z)$, i.e. the sequence
of their zeros and poles that converge to $\zeta_0$.
\begin{thm}\label{t1}
Let  $\zeta_0$ be a zero of one of polynomials $\omega_m(z),
\ \ m=0,1,\allowbreak\ldots,\allowbreak\sigma-1$, and $\zeta_0$ does not
coincide with the poles
$z_{\nu+1},\allowbreak \ldots,z_\lambda$.

If $|\zeta_0|<|z_{\nu+1}|$, then at the point $\zeta_0$ the whole
sequence $\pi^r_{n,\lambda-1}(z)$ converges to $r(z)$.

If
$|\zeta_0|\ge|z_{\nu+1}|$, then $\lim \pi^r_{n,\lambda-1}(\zeta_0)$ does
not exist and also for $|\zeta_0|>|z_{\nu+1}|$
$$
\lim_{n\to\infty} \pi^r_{n,\lambda-1}(\zeta_0)=\infty,\ \ n\in\Lambda_m.
$$
\end{thm}

{\bf Proof.} If $\omega_m(\zeta_0)\ne 0$, then by Theorem~2.3
from~\cite{Adukov} the subsequence
$\left\{\pi^r_{n,\lambda-1}(z)\right\}_{n\in\Lambda_m}$
uniformly converges to $r(z)$ in a neighborhood of $\zeta_0$. It remains
to consider those sequences $\Lambda_m$ for which
$\omega_m(\zeta_0) = 0$. In this case for $n\in\Lambda_m$ we have by
formula (\ref{2})
$$
V_{n+\lambda}(\zeta_0)= D(\zeta_0)\sum_{j=\nu+1}^\lambda
\frac{C_jz_j^{n+\lambda}}{\zeta_0-z_j}.
$$

Taking into account our assumption $|z_{\nu +1}|>|z_{\nu +2}|\ge
\ldots \ge |z_{\lambda}|$, we obtain
\begin{equation}\label{2a}
V_{n+\lambda}(\zeta_0)= z_{\nu+1}^{n+\lambda}D(\zeta_0)\left[
\frac{C_{\nu+1}}{\zeta_0-z_{\nu+1}}+o(1)\right].
\end{equation}
Hence, for all sufficiently large $n\in\Lambda_m$ the polynomial
$z_{\nu+1}^{-(n+\lambda)}V_{n+\lambda}(z)$ is bounded away from zero
in a neighborhood of $\zeta_0$. From the Bezout equation it follows
\begin{equation}\label{3}
r(z)-\pi^r_{n,\lambda-1}(z)=\frac{z^{n+\lambda}}{D(z)V_{n+\lambda}(z)}.
\end{equation}

Thus at $\zeta_0$ we have
$$
r(\zeta_0)-\pi^r_{n,\lambda-1}(\zeta_0)=\frac{\zeta^{n+\lambda}_0}{D^2(\zeta_0)
z_{\nu+1}^{n+\lambda}\left[
\frac{C_{\nu+1}}{\zeta_0-z_{\nu+1}}+o(1)\right]}
$$
Hence, if $|\zeta_0|<|z_{\nu+1}|$, then  at $\zeta_0$ the sequence
$\left\{\pi^r_{n,\lambda-1}(z)\right\}_{n\in\Lambda_m}$  converges
to the function $r(z)$. Therefore the whole sequence
$\pi^r_{n,\lambda-1}(z)$ converges to $r(z)$ at $\zeta_0$.

If $|\zeta_0|\ge|z_{\nu+1}|$, then from (\ref{3}) we obtain for $n\to\infty$,
$n\in\Lambda_m$,
$$
|r(\zeta_0)-\pi^r_{n,\lambda-1}(\zeta_0)|=
\left|\frac{\zeta_0}{z_{\nu+1}}\right|^{n+\lambda}
\left|\frac{1}{D^2(\zeta_)\left[\frac{C_{\nu+1}}{\zeta_0-z_{\nu+1}}+o(1)\right]}
\right|,
$$
i.e. $|r(\zeta_0)-\pi^r_{n,\lambda-1}(\zeta_0)|\to\infty$ for
$|\zeta_0|>|z_{\nu+1}|$ and $|r(\zeta_0)-\pi^r_{n,\lambda-1}(\zeta_0)|\not\to
0$
for $|\zeta_0|=|z_{\nu+1}|$. The theorem is proved. \halmos

Now we construct the examples that illustrate Theorem~\ref{t1}.

\begin{example}\label{nonobstruction}
Let $\lambda=3$, $\nu=2$, $\sigma=2$, $z_1=1, z_2=-1, z_3=1/2$,
$A_1=1$, $A_2=1/18$, $A_3=1$, i.e.
$$
r(z)=\frac{1}{z-1}+\frac{1}{18(z+1)}+\frac{1}{z-1/2}.
$$

Then $C_1=1, C_2=2, C_3=16/9$ and zeros of the polynomials
$\omega_{0,1}(z)$ are $\zeta_0=1/3, \zeta_1=3$. The point
$\zeta_0=1/3$ lies in the unit disk ${\Bbb D}$ and
$|\zeta_0|<|z_3|$. Hence, $\zeta_0$ does not an obstruction for the
poinwise convergence and the sequence $\pi_{n,2}^r(z)$ uniformly
converges to $r(z)$ on compact subsets of the domain ${\Bbb U}_{\Bbb
F}= {\Bbb D}\setminus\{1/3,1/2\}$ and pointwise converges on ${\Bbb
D}\setminus\{1/2\}$.
\end{example}

\begin{example}\label{converg2infty}
Let  $A_1=A_2=A_3=1$,
$$
r(z)=\frac{1}{z-1}+\frac{1}{(z+1)}+\frac{1}{z-1/2}.
$$

Then $C_1=1, C_2=1/9, C_3=16/9$ and $\zeta_0=-4/5$, $\zeta_1=-5/4$.
Now the point $\zeta_0$ is  an obstruction for the poinwise
convergence and
$$
\lim_{m\to\infty}\pi_{2m+1,2}^r(-4/5)=\infty.
$$
 In this case the
sequence $\pi_{n,2}^r(z)$ uniformly converges to $r(z)$
on compact subsets of the domain $\tilde\Bbb U_{\Bbb F}={\Bbb D}\setminus\{1/2,-4/5\}$.
\end{example}

The following example show that an additional limit point can be coincided
with  a nondominant pole.

\begin{example}\label{grmulti}
Let $A_1=2, A_2=2/27, A_3=1$,
$$
r(z)=\frac{2}{z-1}+\frac{2}{27(z+1)}+\frac{1}{z-1/2}.
$$

Then $C_1=1/2, C_2=3/2, C_3=16/9$ and $\zeta_0=1/2, \zeta_1=2$.
The additional limit point $\zeta_0$ coincides with the pole
$z_3$.
\end{example}

The following example show that there are a point $\zeta$ and a sequence
$\Lambda\subset \Bbb N$ such that there is
$\lim_{n\to\infty} \pi_{n,\lambda-1}(\zeta)\ne a(\zeta),
n\in\Lambda$.

\begin{example}\label{noa}
Let $A_1=2/3, A_2=2/9, A_3=1$,
$$
r(z)=\frac{2}{3(z-1)}+\frac{2}{9(z+1)}+\frac{1}{z-1/2}.
$$

Then $C_1=3/2, C_2=1/2, C_3=16/9$ and $\zeta_0=-1/2, \zeta_1=-2$.
By formula~(\ref{3}) we obtain
$\pi_{2m+1,2}^r(-1/2)=0$. Thus $\pi_{2m+1,2}^r(-1/2)$ is a  stationary
sequence and  $\pi_{2m+1,2}^r(-1/2)\ne r(-1/2)=-1$.
\end{example}


\begin{thebibliography}{10}


\bibitem{Adukov}
V.M. Adukov, The uniform convergence of subsequences of the last
intermediate row of the Pad\'e table, {\em J. Approx. Theory} {\bf 122}
(2003), 160--207.

\bibitem{Sidi}
A. Sidi, Quantitative and constructive aspects of the generalized Koenig's
and de Montessus's theorems for Pad\'e approximants, {\em J. Comput.
Appl. Math.} {\bf 29} (1990), 257--291.




\end{thebibliography}
\end{document}